\documentstyle[12pt]{article}
\begin{document}

\author{S.V. Ludkovsky.}

\title{Meta-centralizers of non locally compact group algebras.}

\date{17 April 2013}
\maketitle

\begin{abstract}
Meta-centralizers of non-locally compact group algebras are studied.
Theorems about their representations with the help of families of
generalized measures are proved. Isomorphisms of group algebras are
investigated in relation with meta-centralizers.

 \footnote{key words and phrases: group, algebra,
operator, measure address: Department of Applied Mathematics, \\
Moscow State Technical University MIREA, \\ av. Vernadsky 78, Moscow
119454, Russia\\
sludkowski@mail.ru}

\end{abstract}

\section{Introduction.}
Locally compact group algebras are rather well investigated and play
very important role in mathematics
\cite{fell,hew,nai,kawadamj48,losannm2008}. Left centralizers of
locally compact group algebras were studied in \cite{wendpjm52}. In
all those works Haar measures on locally compact groups were used.
Haar measures are invariant or quasi-invariant relative to left or
right shifts of the entire locally compact group
\cite{bourm70,fell,hew,nai}. According to the A.Weil theorem if a
topological group has a non-trivial borelian measure quasi-invariant
relative to left or right shifts of the entire group, then it is
locally compact. \par On the other hand, the theory of non locally
compact groups and their representations differ drastically from
that of the locally compact case (see
\cite{bana,banasz1,fidal,lujms147:3:08,lujms150:4:08,lurim99} and
references therein). Measures on non locally compact groups
quasi-invariant relative to proper dense subgroups were constructed
in
\cite{beldal,dal,dalshn,lujms147:3:08,lujms150:4:08,lusmldg05,lurim99,lunova2006}.
\par  This article continues investigations of non locally compact group algebras
\cite{bulopnlcga:12,lujms150:4:08,luambp99}. The present paper is
devoted to centralizers of non-locally compact group algebras, which
are substantially different from that of locally compact groups.
Their definition in the non locally compact groups setting is rather
specific and they are already called meta-centralizers. Theorems
about their representations with the help of families of generalized
measures are proved. Isomorphisms of group algebras are investigated
in relation with meta-centralizers. The main results of this paper
are Theorems 8-10 and 14. They are obtained for the first time. \par
Henceforth definitions and notations of \cite{bulopnlcga:12} are
used.
\section{Group algebra}
To avoid misunderstandings we first present our definitions and
notations.
\par {\bf 1. Definition.} Let $\Lambda $ be a directed set and let
$\{ G_{\alpha }: \alpha \in \Lambda \} $ be a family of topological
groups with completely regular (i.e. $T_1\cap T_{3\frac{1}{2}}$)
topologies $\tau _{\alpha }$ satisfying the following restrictions:
\par $(1)$ $\theta ^{\beta }_{\alpha }: G_{\beta }\to G_{\alpha }$
is a continuous algebraic embedding, $\theta ^{\beta }_{\alpha }(
G_{\beta })$ is a proper subgroup in $G_{\alpha }$ for each $\alpha
<\beta \in \Lambda $;
\par $(2)$ $\tau
_{\alpha }\cap {\theta ^{\beta }_{\alpha }(G_{\beta })} \subset
\theta ^{\beta }_{\alpha }(\tau _{\beta })$ and $\theta ^{\beta
}_{\alpha } (G_{\beta })$ is dense in $G_{\alpha }$ for each $\alpha
<\beta \in \Lambda $; then $(\theta ^{\beta }_{\alpha })^{-1}:
\theta ^{\beta }_{\alpha }(G_{\beta }, \tau _{\beta }) \to (G_{\beta
}, \tau _{\beta })$ is considered as the continuous homomorphism;
\par $(3)$ $G_{\alpha } $ is complete relative to the left
uniformity with entourages of the diagonal of the form ${\cal U} =
\{ (h,g): h, g \in G_{\alpha }; h^{-1}g\in U \} $ with neighborhoods
$U$ of the unit element $e_{\alpha }$ in $G_{\alpha }$, $U\in \tau
_{\alpha }$, $e_{\alpha }\in U$;
\par $(4)$ for each $\alpha \in \Lambda $ with $\beta =\phi (\alpha )$
the embedding $\theta ^{\beta }_{\alpha }: (G_{\beta }, \tau _{\beta
})\hookrightarrow (G_{\alpha },\tau _{\alpha })$ is precompact, that
is by our definition for every open set $U$ in $G_{\beta }$
containing the unit element $e_{\beta }$ a neighborhood $V\in \tau
_{\beta }$ of $e_{\beta }$ exists so that $V\subset U$ and $\theta
^{\beta }_{\alpha }(V)$ is precompact in $G_{\alpha }$, i.e. its
closure $cl (\theta ^{\beta }_{\alpha }(V))$ in $G_{\alpha }$ is
compact, where $\phi : \Lambda \to \Lambda $ is an increasing marked
mapping.
\par {\bf 2. Conditions.} Henceforward it is supposed that Conditions $(1-5)$ are
satisfied: \par $(1)$ $\mu _{\alpha }: {\cal B}(G_{\alpha }) \to
[0,1]$ is a probability measure on the Borel $\sigma $-algebra
${\cal B}(G_{\alpha })$ of a group $G_{\alpha }$ from \S 1 with $\mu
_{\alpha }(G_{\alpha })=1$ so that \par $(2)$ $\mu _{\alpha }$ is
quasi-invariant relative to the right and left shifts on $h\in
\theta ^{\beta }_{\alpha }(G_{\beta })$ for each $\alpha <\beta \in
\Lambda $, where $\rho ^{r}_{\mu _{\alpha }}(h,g)=(\mu _{\alpha
}^h)(dg)/\mu (dg)$ and $\rho ^{l}_{\mu _{\alpha }}(h,g)=(\mu
_{\alpha ,h})(dg)/\mu (dg)$ denote quasi-invariance $\mu _{\alpha
}$-integrable factors, $\mu ^h_{\alpha }(S)=\mu (Sh^{-1})$ and $\mu
_{\alpha ,h}(S)=\mu _{\alpha }(h^{-1}S)$ for each Borel subset $S$
in $G_{\alpha }$;
\par $(3)$ a density $\psi _{\alpha }(g) =\mu _{\alpha
}(dg^{-1})/\mu _{\alpha }(dg)$ relative to the inversion exists and
it is $\mu _{\alpha }$-integrable; \par $(4)$ a subset $W_{\alpha
}\in {\cal A}(G_{\alpha })$ exists such that $\rho ^{r}_{\mu
_{\alpha }}(h,g)$ and $\rho ^{l}_{\mu _{\alpha }}(h,g)$ are
continuous on $\theta ^{\beta }_{\alpha }(G_{\beta })\times
W_{\alpha }$ and $\psi _{\alpha }(g)$ is continuous on $W_{\alpha }$
with $\mu _{\alpha }(W_{\alpha })=1$ for each $\alpha \in \Lambda $
with $\beta = \phi (\alpha )$; \par $(5)$ each measure $\mu _{\alpha
}$ is Borel regular and radonian, \\ where the completion of ${\cal
B}(G_{\alpha })$ by all $\mu _{\alpha }$-zero sets is denoted by
${\cal A}(G_{\alpha })$.
\par {\bf 3. Notation.} Denote by $L^1_{G_{\beta }}(G_{\alpha },\mu_{\alpha },{\bf F})$
a subspace in $L^1(G_{\alpha }, \mu _{\alpha },{\bf F})$, which is
the completion of the linear space $L^0(G_{\alpha },{\bf F})$ of all
($\mu _{\alpha }$-measurable) simple functions
$$ f(x)=\sum_{j=1}^n b_j \chi _{B_j}(x),$$ where
$b_j\in {\bf F}$, $~B_j\in {\cal A}(G_{\alpha })$, $~B_j\cap
B_k=\emptyset $ for each $j\ne k$, $ ~ \chi _B$ denotes the
characteristic function of a subset $B$, $~\chi _B(x)=1$ for each
$x\in B$ and $\chi _B(x)=0$ for every $x\in G_{\alpha }\setminus B$,
$~n\in {\bf N}$, where ${\bf F}={\bf R}$ or ${\bf F}={\bf C}$. A
norm on $L^1_{G_{\beta }}(G_{\alpha })$ is by our definition given
by the formula:
$$(1)\quad \| f \|_{L^1_{G_{\beta }}(G_{\alpha })} := \sup_{h\in
\theta ^{\beta }_{\alpha }(G_{\beta })} \| f_h \| _{L^1(G_{\alpha
})} <\infty ,$$ where $f_h(g):=f(h^{-1}g)$ for $h, g \in G_{\alpha
}$, $L^1(G_{\alpha }, \mu _{\alpha },{\bf F})$ is the usual Banach
space of all $\mu _{\alpha }$-measurable functions $u: G_{\alpha
}\to {\bf F}$ such that
$$(2)\quad \| u \| _{L^1(G_{\alpha
})} = \int_{G_{\alpha }} |u(g)|\mu _{\alpha }(dg) <\infty .$$
Suppose that \par $(3)$ $\phi : \Lambda \to \Lambda $ is an
increasing mapping, $\alpha < \phi (\alpha ) $ for each $\alpha \in
\Lambda $. We consider the space \par $(4)$ $L^{\infty
}(L^1_{G_{\beta }}(G_{\alpha },\mu_{\alpha },{\bf F}): \alpha <\beta
\in \Lambda ) := \{ f=(f_{\alpha }: \alpha \in \Lambda ); ~
f_{\alpha }\in L^1_{G_{\beta }}(G_{\alpha },\mu_{\alpha },{\bf F})$
$\mbox{ for each }$ $\alpha \in \Lambda ;$ $\| f \| _{\infty } :=
\sup_{\alpha \in \Lambda } \| f_{\alpha } \|_{L^1_{G_{\beta
}}(G_{\alpha })} <\infty ,$ $\mbox{ where }$ $\beta = \phi (\alpha )
\} .$ \par When measures $\mu _{\alpha }$ are specified, spaces are
denoted shortly by $L^1_{G_{\beta }}(G_{\alpha },{\bf F})$ and
$L^{\infty }(L^1_{G_{\beta }}(G_{\alpha },{\bf F}): \alpha <\beta
\in \Lambda )$ respectively.
\par {\bf 4. Definition.}
Let the algebra ${\cal E} := L^{\infty }(L^1_{G_{\beta }}(G_{\alpha
}, \mu _{\alpha },{\bf F}): \alpha <\beta \in \Lambda )$ be supplied
with the multiplication $f{\tilde \star } u =w$ such that
$$(1)\quad w_{\alpha }(g) = (f_{\beta } {\tilde *} u_{\alpha })(g) =
\int_{G_{\beta }} f_{\beta }(h) u_{\alpha }(\theta ^{\beta }_{\alpha
}(h)g) \mu_{\beta }(dh) $$ for every $f, u \in {\cal E}$ and $g\in
G=\prod_{\alpha \in \Lambda }G_{\alpha }$, where ${\bf F}={\bf R}$
or ${\bf F}={\bf C}$, $\beta =\phi (\alpha )$, $~\alpha \in \Lambda
$.
\par If a bounded linear transformation $T: {\cal E}\to {\cal E}$
satisfies Conditions $(2,3)$:
\par $(2)$ $Tf = (T_{\alpha }f_{\alpha }: \alpha \in \Lambda )$,
$ ~ T_{\alpha }: L^1_{G_{\beta }}(G_{\alpha }, \mu _{\alpha },{\bf
F})\to L^1_{G_{\beta }}(G_{\alpha }, \mu _{\alpha },{\bf F})$ for
each $\alpha \in \Lambda $,
\par $(3)$ $T(f{\tilde \star } u ) = f{\tilde \star } (Tu)$ \\
for each $f, u \in {\cal E}$, then $T$ is called a left
meta-centralizer.
\par {\bf 5. Definitions.} Let $X$ be a topological space, let $C(X,{\bf
R})$ be the space of all continuous functions $f: X\to {\bf R}$,
while $C_b(X,{\bf R})$ be the space of all bounded continuous
functions with the norm
\par $(1)$ $\| f \| := \sup_{x\in X} |f(x)|<\infty $.
\\ Suppose that $\cal F$ is the least $\sigma $-algebra on $X$
containing the algebra $\cal Z$ of all functionally closed subsets
$A=f^{-1}(0)$, $f\in C_b(X,{\bf R})$. A finitely additive
non-negative mapping $m: {\cal F}\to [0,\infty )$ such that
\par $(2)$ $m(A) = \sup \{ m(B): B\in {\cal Z}, B\subset A \} $
\\ for each $A\in \cal F$ is called (a finitely additive) measure.
A generalized measure is the difference of two measures. Denote by
$M(X)=M(X,{\bf R})$ the family of all generalized (finitely
additive) measures.
\par For short "generalized" may be omitted, when $m$ is considered with values
in $\bf R$.
\par {\bf 6. Theorem (A.D. Alexandroff \cite{vardmsb61}).}
{\it $M(X)$ is the topologically dual space to $C_b(X,{\bf R})$,
that is for each bounded linear functional $J$ on $C_b(X,{\bf R})$
there exists a unique generalized (finitely additive) measure $m\in
M(X)$ such that
\par $(1)$ $J(f) = \int_X fdm$ for each $f\in C_b(X,{\bf R})$,
\par each measure $m\in M(X)$ defines a unique continuous linear functional
by Formula $(1)$. Moreover,
\par $(2)$ $\| J \| = \| m \| $.}
\par {\bf 7. Definitions.} A bounded linear functional $J$ on
$C_b(X,{\bf R})$ is called $\sigma $-smooth, if
\par $(1)$ $\lim_n J(f_n)=0$
\\ for each sequence $f_n$ in $C_b(X,{\bf R})$ such that $0\le
f_{n+1}(x)\le f_n(x)$ and $\lim_n f_n(x)=0$ for each point $x\in X$.
The linear space of all $\sigma $-smooth linear functionals is
denoted by $M_{\sigma }(X)=M_{\sigma }(X,{\bf R})$.
\par A bounded linear functional $J$ on $C_b(X,{\bf R})$ is called
tight, if Formula $(1)$ is fulfilled for each net $f_{\alpha }$ in
$C_b(X,{\bf R})$ such that $\| f_{\alpha } \| \le 1$ for each
$\alpha $ and $f_{\alpha }$ tends to zero uniformly on each compact
subset $K$ in $X$. The space of all tight linear functionals is
denoted by $M_t(X)=M_t(X,{\bf R})$.
\par If $m_1, m_2\in M(X)$, then $m=m_1+im_2$ is a complex-valued
measure, their corresponding spaces are denoted by $M(X,{\bf C})$,
$M_{\sigma }(X,{\bf C})=M_{\sigma }(X)+iM_{\sigma }(X)$ and
$M_t(X,{\bf C})= M_t(X)+iM_t(X)$.
\par {\bf 8. Theorem.} {\it Let ${\cal E}$ be a real ${\bf F}={\bf R}$ or complex ${\bf
F}={\bf C}$ algebra (see \S 4), let also $T$ be a left
meta-centralizer on ${\cal E}$. Then there exists a family $\nu =
(\nu _{\alpha }: \alpha \in \Lambda )$ of generalized $\bf F
$-valued measures $\nu _{\alpha }$ on $G_{\alpha }$ of bounded
variation such that
\par $(1)$ $Tf=\nu {\tilde \star } f$, where $$(2)\quad (T_{\alpha } f_{\alpha })(g) = (\nu
_{\beta } {\tilde *} f_{\alpha })(g) = \int_{G_{\beta }} \nu _{\beta
}(dh) f_{\alpha }(\theta ^{\beta }_{\alpha }(h)g)$$ for each $\alpha
\in \Lambda $ and $g\in G_{\alpha }$ with $\beta =\phi (\alpha )$.}
\par {\bf Proof.} For each $\beta \in \Lambda $ and a neutral element
$e_{\beta }\in G_{\beta }$ we consider a basis of its neighborhoods
$ \{ V_{a,\beta }: a\in \Psi _{\beta } \} $ such that $cl_{G_{\alpha
}}\theta ^{\beta }_{\alpha }(V_{a,\beta })$ is compact in
$(G_{\alpha },\tau _{\alpha })$, where $\Psi _{\beta }$ is a set,
$cl_XA$ denotes the closure of a set $A$ in a topological space $X$.
The set $\Psi _{\beta }$ is directed by the inclusion: $a\le b\in
\Psi _{\beta }$ if and only if $V_{b,\beta }\subseteq V_{a,\beta }$.
\par There is a natural continuous linear restriction mapping $p^U_V
: C_b(U,{\bf F})\to C_b(V,{\bf F})$ for each closed subsets $U$ and
$V$ in $G_{\beta }$ such that $V\subset U$, where $p^U_V(f)=f|_V$
for each $f\in C_b(U,{\bf F})$. At the same time, if $U$ is compact,
then each continuous bounded function $g$ on $V$ with values in $\bf
F$ has a continuous extension $\pi ^V_U(g)$ on $U$ with values in
$\bf F$ such that $$ \| g \| _{C_b(V,{\bf F})}\le \| \pi ^V_U(g) \|
_{C_b(U,{\bf F})} \le 2 \| g \| _{C_b(V,{\bf F})}$$ due to
Tietze-Uryson Theorem 2.1.8 \cite{eng}, since $G_{\beta }$ is $T_0$
and hence completely regular by Theorem 8.4 \cite{hew} and each
Huasdorff compact space is normal by Theorems 5.1.1 and 5.1.5
\cite{eng}. Thus there exists a linear continuous embedding $\pi
^V_U: C_b(V,{\bf F})\hookrightarrow C_b(U,{\bf F})$.
\par The probability measure $\mu _{\beta }$ on $G_{\beta }$ is Borel regular
and radonian hence there exists a $\sigma $-compact subset $X_{\beta
}$ in $G_{\beta }$ such that $\mu _{\beta }(X_{\beta })=1$, i.e.
$X_{\beta }$ is the countable union of compact subsets $X_{\beta
,n}$ in $(G_{\beta },\tau _{\beta })$ with $X_{\beta ,n}\subset
X_{\beta ,n+1}$ for each natural number $n$.  \par We put \par $(3)$
$q_{a,\beta } := \chi _{V_{a,\beta }} /\mu _{\beta }(V_{a,\beta })$,
\\ where $\chi _A$ denotes the characteristic function of a subset
$A$ in $G_{\beta }$, $ ~ \chi _A(x)=1$ for each $x\in A$ while $\chi
_A(x)=0$ for each $x\notin A$. In view of Proposition 17.7
\cite{lujms150:4:08} (see also Lemma 13 \cite{bulopnlcga:12}) the
net $ \{ q_{a,\beta }: a\in \Psi _{\beta } \} $ is an approximation
of the identity relative to the convolution: \par $(4)$ $\lim_a
q_{a,\beta } {\tilde
*}f_{\alpha }  =f_{\alpha }$ \\ for each $f_{\alpha }\in L^1_{G_{\beta
}}(G_{\alpha }, \mu _{\alpha },{\bf F})$. From Formulas $(2,4)$ and
4$(1-3)$ it follows that $$(5)\quad T_{\alpha } f_{\alpha } =
T_{\alpha } [\lim_a q_{a,\beta } {\tilde *} f_{\alpha } ] = \lim_a
q_{a,\beta } {\tilde *} [T_{\alpha } f_{\alpha }].$$ Then
$q_{a,\beta } {\tilde *} [T_{\alpha } \cdot ]: L^1_{G_{\beta
}}(G_{\alpha })\to L^1_{G_{\beta }}(G_{\alpha })$ is a continuous
linear operator for each $a\in \Psi _{\beta }$ and $\alpha \in
\Lambda $, particularly, for each $f_{\alpha }$ in the space
$C_b(G_{\alpha },{\bf F})$ of all bounded continuous functions on
$G_{\alpha }$ with values in the field $\bf F$, where \par $(6)$ $
\| f_{\alpha } \| _{C_b} := \sup_{x\in G_{\alpha }} |f_{\alpha
}(x)|<\infty $ \\ for each $f_{\alpha }\in C_b(G_{\alpha },{\bf
F})$. The restriction of each $f_{\alpha }\in C_b(G_{\alpha },{\bf
F})$ on $\theta ^{\beta }_{\alpha }(G_{\beta })$ is bounded and
continuous, while $C_b(G_{\beta },{\bf F})$ is dense in
$L^1_{G_{\gamma }}(G_{\beta },\mu _{\beta },{\bf F})$ with $\gamma =
\phi (\beta )$ (see also Lemma 17.8 and Proposition 17.9
\cite{lujms150:4:08}).
\par This implies that an adjoint operator $B=T^*$ exists
relative to the $\tilde{*}$ multiplication according to the formula:
$$(7)\quad (v_{\beta }{\tilde
*} [T_{\alpha } {\bar f}_{\alpha }])(x) =\int_{G_{\beta }}
v_{\beta }(h)[T_{\alpha } {\bar f}_{\alpha }](\theta^{\beta
}_{\alpha }(h)x)\mu_{\beta }(dh)$$ $$ =: \int_{G_{\beta }} (B_{\beta
} v_{\beta })(h) {\bar f}_{\alpha } (\theta^{\beta }_{\alpha
}(h)x)\mu_{\beta }(dh) $$ for each $v, f \in {\cal E}$, where $x\in
G_{\alpha }$, $~ \bar{z}$ denotes the complex conjugated number of
$z\in \bf C$. The operator $B_{\beta }$ is bounded and linear from
$L^1_{G_{\gamma }}(G_{\beta })$ into itself, since from Formula
$(7)$ the estimate follows:
$$(8)\quad \| B_{\beta } \| \le \sup_{s\in \theta ^{\gamma }_{\beta }(G_{\gamma }),
~ t\in \theta ^{\beta }_{\alpha }(G_{\beta }), ~ 0\ne v_{\beta }\in
L^1_{G_{\gamma }}(G_{\beta }), ~ 0\ne f_{\alpha }\in L^1_{G_{\beta
}}(G_{\alpha }) }$$
$$ \frac{| \int_{G_{\alpha }} \int_{G_{\beta }} v_{\beta } (sh)
[T_{\alpha } {\bar f}_{\alpha }](\theta^{\beta }_{\alpha }(h)tx)
\mu_{\beta }(dh)\mu_{\alpha }(dx)|} {\| v_{\beta } \|
_{L^1_{G_{\gamma }}(G_{\beta })} \| f_{\alpha } \| _{L^1_{G_{\beta
}}(G_{\alpha })} } \le  \| T_{\alpha } \|  <\infty .$$
\par The family of bounded linear operators $ \{ (B_{\beta } q_{a,\beta
}) {\tilde *}: ~ a \in \Psi _{\beta } \} $ from $L^1_{G_{\beta
}}(G_{\alpha })$ into $L^1_{G_{\beta }}(G_{\alpha })$ is pointwise
bounded and hence by the Banach-Steinhaus Theorem (11.6.1)
\cite{nari} it is uniformly bounded: $$(B1)\quad \sup_{a\in \Psi
_{\beta }} \| (B_{\beta }q_{a,\beta } )\tilde{*} \| <\infty .$$
Therefore inequality $(8)$ leads to the conclusion that $B_{\beta }
q_{a,\beta } =: h_{a,\beta }\in L^1_{G_{\gamma }}(G_{\beta }, \mu
_{\beta },{\bf F})$ for every $a\in \Psi _{\beta }$ and $\beta \in
\Lambda $. Each function $h_{a,\beta }$ induces the linear
functional
$$(9)\quad F_{a,\beta } (g_{\beta }) := \int_{G_{\beta }} g_{\beta
}(x) \bar{h}_{a,\beta }(x) \mu _{\beta }(dx).$$
\par Without loss of generality we choose $V_{a,\beta }$ such that
$cl_{G_{\alpha }} V_{a,\beta }$ is compact in $(G_{\alpha }, \tau
_{\alpha })$ for each $a\in \Psi _{\beta }$. Certainly, if $f\in
L^1_{G_{\gamma }}(G_{\beta }, \mu _{\beta },{\bf F})$, then $f\in
L^1(G_{\beta }, \mu _{\beta },{\bf F})$ and $$(10)\quad \| f
\|_{L^1(G_{\beta }, \mu _{\beta },{\bf F})} \le \| f
\|_{L^1_{G_{\gamma }}(G_{\beta }, \mu _{\beta },{\bf F})}<\infty .$$
There is the embedding $C_b(G_{\beta },{\bf F})\subset
L^1_{G_{\gamma }}(G_{\beta }, \mu _{\beta },{\bf F})$ and
$$(11)\quad \| f \|_{L^1_{G_{\gamma }}(G_{\beta }, \mu _{\beta
},{\bf F})}\le \| f \|_{C_b(G_{\beta },{\bf F})} <\infty
$$ for each $f\in C_b(G_{\beta },{\bf F})$, since $\mu _{\beta
}$ is the probability measure on $G_{\beta }$.

\par If $f\in L^1_{G_{\gamma }}(G_{ \beta })$, then
$s\mapsto f{\tilde *}s$ is a continuous linear operator from
$C_b(G_{\beta },{\bf F})$ into $C_b(G_{\beta },{\bf F})$. This
follows from the formulas:
$$(12)\quad (f{\tilde *}s) (g) = \int_{G_{\beta }} f(h) s (hg)\mu _{\beta
}(dh) ,$$ where $g\in G_{\beta }$ and $$\sup_g |(f{\tilde *}s) (g)|
\le \| s \|_{C_b} \int_{G_{\beta }} |f(h) | \mu _{\beta }(dh) \le \|
s \|_{C_b} \| f \| _{L^1(G_{\beta })}\le \| s \|_{C_b} \| f \|
_{L^1_{G_{\gamma }}(G_{\beta })}.$$ It remains to verify that the
function $(f{\tilde *}s) (g)$ is continuous for each $f$ and $s$ as
just above. For the proof consider the term
$$(13)\quad |(f{\tilde *}s) (g_1)- (f{\tilde
*}s) (g_2)| = |\int_{G_{\beta }} f(h) [s (hg_1)-s(hg_2) ] \mu _{\beta
}(dh) |.$$ From $f\in L^1_{G_{\gamma }}(G_{\beta })$ and $s\in
C_b(G_{\beta },{\bf F})$ it follows that for each $\epsilon >0$
there exists a compact subset $V$ in $G_{\beta }$ such that
$\int_{G_{\beta }\setminus V} |f(h)| \mu _{\beta }(dh)  <\epsilon $
and hence $\int_{G_{\beta }\setminus V} |f(h)[s (hg_1)-s(hg_2) ]|
\mu _{\beta }(dh)  <\delta $, where $0<\delta =\epsilon 2 \| s \|
_{C_b}$. Indeed, for each $\delta >0$ there exists a simple function
$q\in L^1_{G_{\gamma }}(G_{\beta })$ such that $ \| f - q
\|_{L^1_{G_{\gamma }}(G_{\beta })}<\delta $ and hence the measure
$|f(h) |\mu _{\beta }(dh)$ is radonian together with $|q(h)|\mu
_{\beta }(dh)$. At the same time, certainly, $\int_{ V} |f(h)| \mu
_{\beta }(dh) \le \| f \|_{L^1(G_{\beta })}$.
\par On the other hand, $[s (hg_1)-s(hg_2) ]$ is uniformly continuous on
$V$ by the variable $h$, since $V$ is compact and $s$ is the
continuous function. For each symmetric open neighborhood $U=U^{-1}$
of the neutral element $e_{\beta }$ in $G_{\beta }$ there exists a
finite family of elements $p_1,...,p_n\in G_{\beta }$ such that
$V\subset p_1U\cup ... \cup p_nU$, since $V$ is compact. Thus
$VU\subset p_1U^2\cup ... \cup p_nU^2$. Consider a family of
symmetric open neighborhoods $U_k=U_k^{-1}$ of $e_{\beta }$ such
that $\{ p_kU_k: k\in \omega \} $ is a covering of $V$ and $|s
(hg_1)-s(hg_2) |<\epsilon $ for each $h\in p_kU_k$ and $g_1, g_2\in
U_k$, where $p_k\in G_{\beta }$ for each $k$, whilst $\omega $ is an
ordinal. The covering $p_kU_k$ of $V$ has a finite subcovering for
$k\in M$, where $M$ is a finite subset in $\omega $. Thus for each
$\epsilon
>0$ there exists a symmetric neighborhood $U \subseteq \bigcap_{k\in
M} U_k$ of $e_{\beta }$ such that $|s (hg_1)-s(hg_2) |<\epsilon $
for each $h\in V$ and $g_1, g_2\in U$. Therefore, $$|(f{\tilde *}s)
(g_1)- (f{\tilde *}s) (g_2)| \le \delta + \epsilon \| f \| _{L^1} =
\epsilon (\| f \| _{L^1}+ 2 \| s \| _{C_b})$$ for each $g_1, g_2\in
U$. Thus \par $(14)$ $f{\tilde *}s \in C_b(G_{\beta },{\bf F})$ for
each $f\in L^1_{G_{\gamma }}(G_{\beta },{\bf F})$ and $s\in
C_b(G_{\beta },{\bf F})$.
\par This implies that \par $(15)$ $C_b(G_{\beta },{\bf F})\ni s\mapsto (f{\tilde
*}s)(e_{\beta }) \in \bf F$ \\ is the continuous linear functional on
$C_b(G_{\beta },{\bf F})$. In particular each operator $(B_{\beta }
q_{a,\beta }) {\tilde *}$ indices the continuous linear functional
\par $(16)$ $J_{a, \beta }(s) = [(B_{\beta } q_{a,\beta }) {\tilde
*}s](e_{\beta })$ on $C_b(G_{\beta },{\bf F})$.
\par There are the inclusions $M_t(X)\subset M_{\sigma }(X)\subset M(X)$
(see \S I.4 \cite{vardmsb61} and Definitions 5, 7 and Theorem 6
above) and for $X=G_{\beta }$ in particular. On the other hand, each
$w_{a,\beta }(dx) := (B_{\beta } q_{a,\beta })(x) \mu _{\beta }(dx)$
is the radonian measure on $G_{\beta }$, i.e. belongs to the space
$M_t(G_{\beta },{\bf F})$ of radonian measures on $G_{\beta }$.
\par Let $\Phi _{\beta }$ be a family of all left-invariant
pseudo-metrics on $(G_{\beta },\tau _{\beta })$ providing its left
uniformity denoted by ${\cal L}_{\beta }$ (see \S 8.1.7 \cite{eng}
and Condition 1$(3)$). This means that each $\kappa \in \Phi _{\beta
}$ satisfies the restrictions: \par $(P1)$ $\kappa (x,y)\ge 0$, \par
$(P2)$ $\kappa (x,x)=0$,
\par $(P3)$ $\kappa (x,y)=\kappa (y,x)$,
\par $(P4)$ $\kappa (x,y)\le \kappa (x,z)+\kappa (z,y)$
\par $(P5)$ $\kappa (zx,zy) = \kappa (x,y)$ for each $x, y, z \in G_{\beta
}$. \par The family $\Phi _{\beta }$ is directed: $\kappa _1\le
\kappa \in \Phi _{\beta }$ if and only if $\kappa _1(x,y)\le \kappa
(x,y)$ for each $x, y \in G_{\beta }$; without loss of generality
for each $\kappa , \kappa _1\in \Phi _{\beta }$ there exists $\kappa
_2\in \Phi _{\beta }$ such that $\kappa \le \kappa _2$ and $\kappa
_1\le \kappa _2$, since $\kappa + \kappa _1\in \Phi _{\beta }$. Each
pseudo-metric $\kappa \in \Phi _{\beta }$ defines the equivalence
relation: $x\Xi _{\kappa }y$ if and only if $\kappa (x,y)=0$. Then
as the uniform space $(G_{\beta }, {\cal L}_{\beta })$ has the
projective limit decomposition (i.e. the limit of the inverse
mapping system) $$G_{\beta } = \lim \{ G_{\beta ,\kappa }, \pi
^{\kappa }_{\omega }, \Phi _{\beta } \} ,$$ where $G_{\beta , \kappa
} := G_{\beta }/\Xi _{\kappa }$ denotes the quotient uniform space
with the quotient uniformly, $\pi _{\kappa }$ is a uniformly
continuous mapping from $G_{\beta }$ onto $G_{\beta , \kappa }$, $ ~
\pi ^{\kappa }_{\omega }$ are uniformly continuous mappings from
$G_{\beta , \kappa }$ onto $G_{\beta , \omega } $ for each $\omega
\le \kappa \in \Psi _{\beta }$ such that $\pi ^{\omega }_{\xi }\circ
\pi ^{\kappa }_{\omega } = \pi ^{\kappa }_{\xi }$ and $\pi _{\omega
} = \pi ^{\kappa }_{\omega }\circ \pi _{\kappa }$ for each $\xi \le
\omega \le \kappa \in \Phi _{\beta }$ (see \S \S 8.2.B, 2.5.F and
Proposition 2.4.2 \cite{eng} or \cite{isb5}). Moreover, the equality
is satisfied: $ \{ y\in G_{\beta }: ~ x\Xi _{\kappa } y \} = x\Omega
_{\beta , \kappa }$ with $\Omega _{\beta , \kappa } := \{ y\in
G_{\beta }: ~ e_{\beta }\Xi _{\kappa } y \} $, since $\kappa
(x,y)=0$ if and only if $\kappa (e_{\beta },x^{-1}y)=0$ by Property
$(P5)$, where $e_{\beta }$ denotes the neutral element in the group
$G_{\beta }$. That is, $G_{\beta, \kappa }$ is called the
homogeneous quotient uniform space.
\par At the same time the $\sigma $-compact subset $X_{\beta }$ is
dense in $G_{\beta }$, since $\mu _{\beta }(U)>0$ for each open
subset $U$ in $G_{\beta }$, but $\mu _{\beta }(X_{\beta }) = \mu
_{\beta }(G_{\beta })=1$ (see the proof above). Therefore, $\pi
_{\kappa }(X_{\beta })$ is dense in $G_{\beta ,\kappa }$. Then $\pi
_{\kappa }(X_{\beta ,n})$ is compact for each $\kappa \in \Phi
_{\beta }$ as the continuous image of the compact space according to
Theorem 3.1.10 \cite{eng}, consequently, $\pi _{\kappa }(X_{\beta
})= \bigcup_{n=1}^{\infty } \pi _{\kappa }(X_{\beta ,n})$ is $\sigma
$-compact. On the other hand, $G_{\beta ,\kappa }$ is metrizable and
complete, since $(G_{\beta }, {\cal L}_{\beta })$ is complete.
Therefore, the topological space $\pi _{\kappa }(X_{\beta })$ is
separable, since each $\pi _{\kappa }(X_{\beta ,n})$ is separable by
Theorems 4.3.5 and 4.3.27 \cite{eng} and $\pi _{\kappa }(X_{\beta
})= \bigcup_{n=1}^{\infty } \pi _{\kappa }(X_{\beta ,n})$. This
implies that each metrizable space $G_{\beta ,\kappa }$ is separable
and complete. \par The spaces $C_b(G_{\beta },{\bf F})$ and
$C_b^*(G_{\beta },{\bf F})$ form the dual pair (see \S \S 9.1 and
9.2 \cite{nari}). Then we get that the space of bounded continuous
functions $C_b(G_{\beta },{\bf F})$ has the inductive limit
representation $C_b(G_{\beta },{\bf F})= ind-\lim_{\Phi _{\beta }}
C_b(G_{\beta ,\kappa },{\bf F})$, while  its topologically dual
space has the projective limit decomposition $C_b^*(G_{\beta },{\bf
F})= pr-\lim_{\Phi _{\beta }} C_b^*(G_{\beta ,\kappa },{\bf F})$
(see \S \S 9.4, 9.9, 12.2, 12.202 \cite{nari} and also the note
after Theorem 2.5.14 in \cite{eng}). This implies that $\nu _{\beta
}\in M(G_{\beta },{\bf F})$ if and only if
$$(M1)\quad \nu _{\beta } = \lim \{ \nu _{\beta ,\kappa }, \pi
^{\kappa }_{\omega }, \Phi _{\beta } \} ,$$ where $\nu _{\beta
,\kappa }\in M(G_{\beta ,\kappa },{\bf F})$ for each $\kappa \in
\Phi _{\beta }$ so that
$$(M2)\quad \nu _{\beta }(\pi _{\omega }^{-1}(C))=\nu _{\beta , \omega
}(C)\mbox{ and } \nu _{\beta , \kappa }((\pi ^{\kappa }_{\omega
})^{-1}(C))=\nu _{\beta , \omega }(C)$$ for every $C\in {\cal
B}(G_{\beta ,\omega })$ and $\omega \le \kappa \in \Phi _{\beta }$.
\par Then we consider the measure net $\{ w_{a, \beta , \kappa }: a \in \Psi _{\beta } \} $
for each $\kappa \in \Phi _{\beta }$ corresponding to measures
$w_{a,\beta }(dx) = (B_{\beta }q_{a,\beta })(x) \mu _{\beta }(dx)$
according to Formula $(M2)$, where $x\in G_{\beta }$. Since the
measure $w_{a,\beta }(dx)$ is absolutely continuous relative to the
radonian measure $\mu _{\beta }$, then $w_{a,\beta }$ is also
radonian. Therefore, there is the inclusion $\{ w_{a, \beta , \kappa
}: a \in \Psi _{\beta } \} \subset M_t(G_{\beta ,\kappa },{\bf F})$
and it is known that $M_t(Y,{\bf F})\subset M_{\sigma }(Y,{\bf
F})\subset M(Y,{\bf F})$ for a completely regular topological space
$Y$. Thus the measure net $\{ w_{a, \beta }: a \in \Psi _{\beta }
\}$ weakly converges to some measure $\nu _{\beta }$ in $M(G_{\beta
},{\bf F})$ if and only if the net $\{ w_{a, \beta , \kappa }: a \in
\Psi _{\beta } \} $ weakly converges in $M(G_{\beta ,\kappa },{\bf
F})$ for each $\kappa \in \Phi _{\beta }$ according to Theorem 2.5.6
and Corollary 2.5.7 \cite{eng}. The net $\{ w_{a, \beta }: a \in
\Psi _{\beta } \} $ is norm bounded, since
$$\| B_{\beta }q_{a,\beta } \|
_{L^1(G_{\beta })} \le \sup \{ \| (B_{\beta }q_{a,\beta })\tilde{*}
f_{\alpha } \| _{L^1_{G_{\beta }}(G_{\alpha })}: ~ f_{\alpha } \in
L^1_{G_{\beta }}(G_{\alpha }), ~ \| f_{\alpha } \|_{L^1_{G_{\beta
}}(G_{\alpha })}\le 1 \} $$
$$ = \sup \{ \| q_{a,\beta }\tilde{*}
(T_{\alpha } f_{\alpha }) \| _{L^1_{G_{\beta }}(G_{\alpha })}: ~
f_{\alpha } \in L^1_{G_{\beta }}(G_{\alpha }), ~ \| f_{\alpha }
\|_{L^1_{G_{\beta }}(G_{\alpha })}\le 1 \} \le $$
$$ \| T_{\alpha } \|  \sup \{ \| q_{a,\beta }\tilde{*}
g_{\alpha } \| _{L^1_{G_{\beta }}(G_{\alpha })}: ~ g_{\alpha } \in
L^1_{G_{\beta }}(G_{\alpha }), ~ \| g_{\alpha } \|_{L^1_{G_{\beta
}}(G_{\alpha })}\le 1 \}$$ $$ \le \| T_{\alpha } \| <\infty  ,
\mbox{ since }$$
$$ \| u_{\beta } \tilde{*}g_{\alpha } \|_{L^1_{G_{\beta }}(G_{\alpha
})} \le \| u \| _{L^1(G_{\beta })} \| g_{\alpha } \|_{L^1_{G_{\beta
}}(G_{\alpha })} $$ for each $u\in L^1(G_{\beta })$ and $g_{\alpha }
\in L^1_{G_{\beta }}(G_{\alpha })$ (see Lemma 17.2
\cite{lujms150:4:08}). This implies that for each $\epsilon
>0$ and $\kappa \in \Phi _{\beta }$ there exists a compact set
$K_{\epsilon ,\kappa }$ in $G_{\beta ,\kappa }$ such that $w_{a,
\beta , \kappa }(G_{\beta ,\kappa }\setminus K_{\epsilon ,\kappa
})<\epsilon $ for each $a\in \Psi _{\beta }$, since $\mu _{\beta
,\kappa }$ as the image of $\mu _{\beta }$ is the radonian measure
on the complete separable metric space $G_{\beta ,\kappa }$ and each
measure $w_{a, \beta , \kappa }$ is absolutely continuous relative
to $\mu _{\beta ,\kappa }$ (see also Theorem 1.2 \cite{dal} and
Formulas $(M1,M2)$).
\par Applying Theorems either II.24 and II.27 or II.30 \cite{vardmsb61}
we get that a measure $\nu _{\beta ,\kappa }\in M_{\sigma }(G_{\beta
,\kappa },{\bf F})$ exists such that the net $w_{a,\beta ,\kappa }$
weakly converges to $\nu _{\beta ,\kappa }$ for each $\beta \in
\Lambda $ and $\kappa \in \Phi _{\beta }$. Thus using Formulas
$(M1,M2)$ we have deduced that
$$(17)\quad \lim_a J_{a, \beta }(f) = \int_{G_{\beta }} fd\nu _{\beta }$$
for each $f\in C_b(G_{\beta },{\bf F})$. The variation of $\nu
_{\beta }$ is finite and $M(G_{\beta },{\bf F})$ is the Banach space
relative to the variation norm according to Theorems I.2 and I.3
\cite{vardmsb61}.

\par Let $x\in C_b(G_{\beta },{\bf F})$ and $y\in C_b(G_{\gamma },{\bf F})$,
we consider the function $$(18)\quad z(g)=\int_{G_{\gamma }}
y(h)x(\theta ^{\gamma }_{\beta }(h)g)\mu _{\gamma }(dh).$$ It
evidently exists and is $\mu _{\beta }$-measurable, since $\mu
_{\gamma }(G_{\gamma })=1$, consequently, $$\sup_{g\in G_{\beta }}
|\int_{G_{\gamma }} y(h)x(\theta ^{\gamma }_{\beta }(h)g)\mu
_{\gamma }(dh)|\le \| y \|_{C_b(G_{\gamma },{\bf F})} \| x
\|_{C_b(G_{\beta },{\bf F})}.$$ Moreover, $z\in C_b(G_{\beta },{\bf
F})\subset L^1_{G_{\gamma }}(G_{\beta })$ due to the latter
inequality and Properties $(11,14)$ (see above). Since $\nu _{\beta
}$ is the weak limit of the net $J_{a, \beta }$, then for each
$\epsilon >0$ there exists $b\in \Psi _{\beta }$ such that
$$(19)\quad |\int_{G_{\beta }} z(g)\nu _{\beta } (dg) - \int_{G_{\beta }} z(g)
(B_{\beta } q_{a,\beta })(g)\mu _{\beta } (dg)|<\epsilon $$ for each
$a>b$. In view of the Fubini theorem the latter inequality implies
that
$$(20)\quad |\int_{G_{\gamma }} y(h)\mu _{\gamma }(dh) \int_{G_{\beta }} x(\theta ^{\gamma }_{\beta }(h)g)
\nu _{\beta } (dg) -$$ $$ \int_{G_{\gamma }} y(h)\mu _{\gamma }(dh)
\int_{G_{\beta }} x(\theta ^{\gamma }_{\beta }(h)g) (B_{\beta }
q_{a,\beta })(g) \mu _{\beta } (dg)|\le \epsilon $$ for each $a>b$.
Therefore, $T_{\alpha } x(g) = (\nu _{\beta } {\tilde *} x)(g)$ for
each $x\in C_b(G_{\beta },{\bf F})\cap [(\theta ^{\beta }_{\alpha
})^{-1}(C_b(G_{\alpha },{\bf F}))]$ and $g\in G_{\beta }$. If
$f_{\alpha }\in C_b(G_{\alpha },{\bf F})$, then its restriction
$f_{\alpha }|_{\theta ^{\beta }_{\alpha }(G_{\beta })}$ is
continuous and bounded, that is $f_{\alpha }\circ (\theta ^{\beta
}_{\alpha })^{-1}$ is continuous and bounded on $(G_{\beta }, \tau
_{\beta })$ due to 1$(2)$. Moreover, the function $\psi _g(h) :=
f_{\alpha }(\theta ^{\beta }_{\alpha }(h)g)$ is continuous and
bounded by $h\in G_{\beta }$ for each $g\in G_{\alpha }$. Hence
$$(21)\quad (\nu _{\beta } {\tilde *} \psi _g)(s) =\int_{G_{\beta
}}f_{\alpha }(\theta ^{\beta }_{\alpha }(hs)g) \nu _{\beta } (dh) =
[\nu _{\beta }{\tilde
*} f_{\alpha }](\theta ^{\beta }_{\alpha }(s)g)$$ is defined for
each $s\in G_{\beta }$ and $g\in G_{\alpha }$, particularly for
$s=e_{\beta }$.
\par By the conditions of this theorem $T_{\alpha }: L^1_{G_{\beta
}}(G_{\alpha })\to L^1_{G_{\beta }}(G_{\alpha })$ is the continuous
linear operator. There is also the inclusion $C_b(G_{\alpha },{\bf
F})\subset L^1_{G_{\beta }}(G_{\alpha },\mu_{\alpha },{\bf F})$ so
that $C_b(G_{\alpha },{\bf F})$ is dense in $L^1_{G_{\beta
}}(G_{\alpha },\mu_{\alpha },{\bf F})$, since $\mu _{\alpha
}(X_{\alpha })=\mu _{\alpha }(G_{\alpha })=1$ with the $\sigma
$-compact subset $X_{\alpha }$ in $G_{\alpha }$ (see also Lemma 17.8
and Proposition 17.9 \cite{lujms150:4:08} and Property $(14)$
above). Let $f_{\alpha }\in L^1_{G_{\beta }}(G_{\alpha },\mu_{\alpha
},{\bf F})$ and we take any sequence of bounded continuous functions
$f_{\alpha ,n}\in C_b(G_{\alpha },{\bf F})$ converging to $f_{\alpha
}$ in $L^1_{G_{\beta }}(G_{\alpha },\mu_{\alpha },{\bf F})$. We have
$$(22)\quad \lim_a (B_{\beta } q_{a,\beta }) {\tilde *}f_{\alpha ,n}=
f_{\alpha } \mbox{  and  }\lim_n f_{\alpha ,n} = f_{\alpha }$$ in
$L^1_{G_{\beta }}(G_{\alpha },\mu_{\alpha },{\bf F})$.  Then
$$(23)\quad \| (B_{\beta } q_{a,\beta }) {\tilde
*}f_{\alpha ,n} - (B_{\beta } q_{b,\beta }) {\tilde *}f_{\alpha ,m}
\|_{L^1_{G_{\beta }}(G_{\alpha })} \le $$  $$\| (B_{\beta }
q_{a,\beta } - B_{\beta } q_{b,\beta }) {\tilde
*} f_{\alpha ,n} \|_{L^1_{G_{\beta }}(G_{\alpha })} +
\| (B_{\beta } q_{b,\beta })\tilde{*} \| \| f_{\alpha ,n} -
f_{\alpha ,m} \|_{L^1_{G_{\beta }}(G_{\alpha })} ,$$ consequently,
for each $\epsilon >0$ there exist $a_0\in \Psi _{\beta }$ and
$n_0\in {\bf N}$ such that
$$(24)\quad \| (B_{\beta } q_{a,\beta }) {\tilde *}f_{\alpha ,n} - (B_{\beta }
q_{b,\beta }) {\tilde *}f_{\alpha ,m} \|_{L^1_{G_{\beta }}(G_{\alpha
})} <\epsilon $$ for each $a, b>a_0$ and $n, m>n_0$ (see Lemma 17.2
and Proposition 17.7 \cite{lujms150:4:08} and Formula $(B1)$ above).
That is the net $ \{ (B_{\beta } q_{a,\beta }) {\tilde
*} f_{\alpha ,n}: (a,n) \} $ is fundamental (i.e. of the Cauchy
type) in the Banach space $L^1_{G_{\beta }}(G_{\alpha })$, where
$(a,n)\le (b,m)$ if $a\le b$ and $n\le m$. Therefore the limit
exists
$$(25)\quad T_{\alpha }f_{\alpha } = \lim_{a,n} (B_{\beta } q_{a,\beta
}) {\tilde
*}f_{\alpha ,n}=\lim_n \lim_a (B_{\beta } q_{a,\beta }) {\tilde
*}f_{\alpha ,n}= \lim_n \nu_{\beta }{\tilde *} f_{\alpha ,n} = \nu_{\beta }{\tilde *} f_{\alpha }.$$
Thus $$T_{\alpha } f_{\alpha } = \nu _{\beta } {\tilde
*} f_{\alpha }$$ for each $f_{\alpha }\in L^1_{G_{\beta }}(G_{\alpha })$
as well, that is, Formulas $(1,2)$ are fulfilled.
\par {\bf 9. Theorem.} {\it Let suppositions of Theorem 8 be satisfies.
Then the statement of Theorem 8 is equivalent to the following:
\par $(1)$ relative to the strong operator topology the set of all
convolution operators of the form 8$(1,2)$ on ${\cal E}:= L^{\infty
}(L^1_{G_{\beta }}(G_{\alpha }, \mu _{\alpha },{\bf F}): \alpha
<\beta \in \Lambda )$ with values in $\cal E$ is a closed subset of
the ring of all bounded linear operators from ${\cal E}$ into ${\cal
E}$.}
\par {\bf Proof.} $(8\Rightarrow 9)$. Let $\nu _{a,\beta }{\tilde *}$ be a net of convolution
operators converging to an operator $T_{\alpha }: L^1_{G_{\beta
}}(G_{\alpha })\to L^1_{G_{\beta }}(G_{\alpha })$ in the strong
operator topology for each $\alpha \in \Lambda $, hence $T$ is the
left meta-centralizer on ${\cal E}$, since each operator $\{ \nu
_{a,\beta }\tilde{*}: ~ \alpha \in \Lambda  , \beta = \phi (\alpha )
\} $ is the left meta-centralizer.
\par $(9\Rightarrow 8)$. From the proof of Theorem 8 we analogously
get
$$T_{\alpha } f_{\alpha } = \lim_a \nu _{a,\beta } {\tilde *}
f_{\alpha }$$ for each $\alpha \in \Lambda $ and $f_{\alpha }\in
L^1_{G_{\beta }}(G_{\alpha }, \mu _{\alpha },{\bf F})$ with $\beta =
\phi (\alpha )$, where $\nu _{a,\beta }\in M(G_{\beta },{\bf F})$
for each $\beta \in \Lambda $ and $a\in \Psi _{\beta }$
consequently, $T= (T_{\alpha }: \alpha )$ is the convolution
operator.

\par {\bf 10. Theorem.} {\it Let $S$ be a bounded linear mapping of ${\cal E}$ (see \S 4)
into itself such that $Sf= (S_{\alpha }f_{\alpha }: \alpha \in
\Lambda )$ with $S_{\alpha }: L^1_{G_{\beta }}(G_{\alpha })\to
L^1_{G_{\beta }}(G_{\alpha })$ for each $\alpha \in \Lambda $ with
$\beta = \phi (\alpha )$. Then the following statements $(i)$ and
$(ii)$ are equivalent: \par $(i)$ an operator $S$ has the form
\par $(1)$ $S=p{\hat U}_a$ for some marked elements $a\in G
:= \prod_{\alpha \in \Lambda }G_{\alpha }$ and $p= \{ p_{\alpha }: |p_{\alpha }|=1 ~
\forall \alpha \in \Lambda \} \in {\bf F}^{\Lambda }$, that is
\par $(2)$ $S_{\alpha }f_{\alpha }(x) = p_{\alpha }{\hat U}_{a_{\beta
}}f_{\alpha }(x)$ for any $\alpha \in \Lambda $ with $\beta =\phi
(\alpha )$ and each $x\in G_{\alpha }$, where \par $(3)$ ${\hat
U}_{g_{\beta }} f_{\alpha }(x) = f_{\alpha }(\theta ^{\beta
}_{\alpha }(g_{\beta })x)$ for each $g_{\beta }\in G_{\beta }$ and $
x\in G_{\alpha }$;
\par $(ii)$ (4) $S$ is a left meta-centralizer and \par $(5)$ $ \| S_{\alpha } f_{\alpha } \| =
\| f_{\alpha } \| $ for every $f_{\alpha }\in L^1_{G_{\beta
}}(G_{\alpha })$ and  $\alpha \in \Lambda $ with $\beta = \phi
(\alpha )$.}
\par {\bf Proof.} The $\bf F$-linear span of the set of all non-negative
functions $f\in L^1_{G_{\beta }}(G_{\alpha },\mu_{\alpha },{\bf F})$
is dense in $L^1_{G_{\beta }}(G_{\alpha },\mu_{\alpha },{\bf F})$.
Therefore, each bounded linear operator $S_{\alpha }$ can be written
in the form $S_{\alpha } = S_{1,\alpha }+iS_{2,\alpha
}=S^+_{1,\alpha }-S^-_{1,\alpha }+iS^+_{2,\alpha }-iS^-_{2,\alpha
}$, where $S^+_{k,\alpha }f\ge 0$ and $S^-_{k,\alpha }f\ge 0$ for
$k=1,2$ and each $f\in P_{\alpha }$, $ ~ S_{k,\alpha }=
S^+_{k,\alpha }-S^-_{k,\alpha }$, where $~ P_{\alpha }$ denotes the
cone of functions in $L^1_{G_{\beta }}(G_{\alpha },\mu_{\alpha
},{\bf F})$ non-negative $\mu _{\alpha }$-almost everywhere on
$G_{\alpha }$. Certainly over the real field additives $S_{2,\alpha
}^{\pm }$ vanish. In view of Theorem 11 \cite{bulopnlcga:12} there
exist $a^+_k\in G$ and $p^+_k= \{ p^+_{k,\alpha }: p^+_{k,\alpha }>0
~ \forall \alpha \in \Lambda \} \in {\bf R}^{\Lambda }$ such that
$S^+_{k,\alpha }f_{\alpha }(x) = p^+_{k,\alpha }{\hat
U}_{a^+_{k,\beta }}f_{\alpha }(x)$ and analogously for
$S^-_{k,\alpha }$ for each $k=1, 2$. \par Suppose that $a^t_k\ne
a^s_l$ for some $t, s \in \{ +, - \}$ and $k, l\in \{ 1, 2 \} $,
then there exists $\alpha \in \Lambda $ such that $a^t_{k,\beta }\ne
a^s_{l,\beta }$ with $\beta =\phi (\alpha )$. On the other hand, we
have $S_{k,\alpha } f_{\alpha } = S^+_{k,\alpha }f_{\alpha } -
S^-_{k,\alpha }f_{\alpha } = p^+_{k,\alpha }f_{\alpha }(\theta
^{\beta }_{\alpha }(a^+_{k,\beta })x) - p^-_{k,\alpha }f_{\alpha
}(\theta ^{\beta }_{\alpha }(a^-_{k,\beta })x)$ for each $f_{\alpha
} \in L^1_{G_{\beta }}(G_{\alpha },\mu _{\alpha },{\bf F})$, since
$f_{\alpha }=[f^+_{1,\alpha } - f^-_{1\alpha }] + i [f^+_{2,\alpha }
- f^-_{2,\alpha }]$ , where $f^+_{k,\alpha }(x) = \max(f_{k,\alpha
}(x),0)$ for every $k=1, 2$ and $x\in G_{\alpha }$, $ ~
f^+_{k,\alpha }, f^-_{k,\alpha } \in P_{\alpha }$. Then if $U$ is an
open subset in $G_{\alpha }$ such that $\theta ^{\beta }_{\alpha
}(a^s_{k,\beta })U\cap \theta ^{\beta }_{\alpha }(a^t_{l,\beta })U =
\emptyset $ for every $k, l=1,2$ and $t, s \in \{ +, - \} $, then  $
\| S_{\alpha }\chi _U \| = \sum_{k=1}^2\sum_{t\in \{ +, - \} }
(|p^t_{k,\alpha }| \| {\hat U}_{a^t_{k,\beta }}\chi _U \|)$. If the
interior of the intersection $ \cap_{k=1}^2\cap_{t\in \{ +, - \} }
(\theta ^{\beta }_{\alpha }(a^t_{k,\beta })U)$ is non-void, then $
\| S_{\alpha }\chi _U \| < \sum_{k=1}^2\sum_{t\in \{ +, - \} }
(|p^t_{k,\alpha }| \| {\hat U}_{a^t_{k,\beta }}\chi _U \|)$, since
$\mu _{\alpha }(V)>0$ for each open subset $V$ in $G_{\alpha }$,
consequently, $S_{\alpha }$ is not an isometry.
\par Therefore, if $S$ satisfies Conditions $ii(4,5)$, then $a^t_{k,\beta }=a^s_{l,\beta }$
for each $t, s\in \{ +, - \} $ and $k, l \in \{ 1, 2 \} $. Thus
$(S_{\alpha }f_{\alpha }) = p_{\alpha }{\hat U}_{a_{\beta
}}f_{\alpha }(x)$ for any $\alpha \in \Lambda $ and each $x\in
G_{\alpha }$, where $p_{\alpha } = p^+_{1,\alpha } - p^-_{1,\alpha }
+ ip^+_{2,\alpha } - ip^-_{2,\alpha } $. Naturally, in the case
${\bf F}={\bf R}$ the terms $p_2^{\pm }$ vanish. In view of Lemma 7
\cite{bulopnlcga:12} ${\hat U}_a$ is the isometry. Since $S$
preserves norms, then $|p_{\alpha }|=1$ for each $\alpha $. \par
Vice versa Conditions $i(1-3)$ imply $ii(4,5)$ due to Lemma 7
\cite{bulopnlcga:12}.

\par {\bf 11. Lemma.} {\it Let ${\hat U}_c$ be a left translation on $\cal E$
as in \S 10, let also $T: {\cal E}\to {\cal F}$ be an isomorphism of
normed algebras such that $Tf=(T_{\alpha }f_{\alpha }: \alpha \in
\Lambda )$, $T_{\alpha }: L^1_{G_{\beta }}(G_{\alpha }, \mu _{\alpha
},{\bf F})\to L^1_{H_{\beta }}(H_{\alpha }, \lambda _{\alpha },{\bf
F})$ and $\| T_{\alpha } \| \le 1$ for each $\alpha $, where ${\cal
F} = L^{\infty }(L^1_{H_{\beta }}(H_{\alpha }, \lambda _{\alpha
},{\bf F}): \alpha <\beta \in \Lambda )$. If ${\hat K}_c = T {\hat
U}_cT^{-1}$, then there exist mappings of groups $\xi : G\to H$ and
$p: G\to {\bf F}^{\Lambda }$ such that
\par $(1)$ ${\hat K}_c=p_c{\hat V}_t$ for $t=\xi (c)$ and \par $(2)$ $p_c = \{
p_{c,\alpha }: |p_{c,\alpha }|=1 ~ \forall \alpha \in \Lambda \} \in
{\bf F}^{\Lambda }$, where ${\hat V}_d$ denotes the left translation
operator on $\cal F$, $ ~ c\in G$.}
\par {\bf Proof.} We have $T(f{\tilde \star } u) = (Tf){\tilde \star }
(Tu)$ for each $u, f \in {\cal E}$ and $T^{-1}(g{\tilde \star } v) =
(T^{-1}g){\tilde \star } (T^{-1}v)$ for each $v, g \in {\cal F}$.
One can take the approximate identity $\{ q_{a,\beta }: a\in \Psi
_{\beta } \} $ as in \S 8 and consider functions $s_{a,\beta } =
T_{\beta }q_{a,\beta }$. The operator $T$ is bijective and
continuous from ${\cal E}$ onto ${\cal F}$, where ${\cal E}$ and
${\cal F}$ as linear normed spaces are complete. According to the
Banach theorem IV.5.4.3 \cite{kolmfom} (or see \cite{danschw}) the
inverse operator $T^{-1}$ is also bounded. Due to Formulas 8$(7,8)$
there exists the adjoint operator $({\hat K}_{c_{\gamma }})^*$
relative to the $\tilde{*}$ multiplication for each $c\in G$ and
$\gamma \in \Lambda $. For each $f, g\in {\cal F}$, $\gamma =\phi
(\beta )$ and $\beta = \phi (\alpha )$ the limit exists
$$({\hat K}_{c_{\gamma }}f_{\beta })\tilde{*} g_{\alpha
} =f_{\beta }\tilde{*} [({\hat K}_{c_{\gamma }})^* g_{\alpha }] =
\lim_a f_{\beta }\tilde{*} \{ s_{a,\beta }\tilde{*}[({\hat
K}_{c_{\gamma }})^* g_{\alpha }] \} $$
$$ = f_{\beta }\tilde{*} \{ \lim_a ({\hat K}_{c_{\gamma }} s_{a,\beta })\tilde{*} g_{\alpha } \} =
f_{\beta }\tilde{*} \{ \lim_a (T_{\beta } {\hat U}_{c_{\gamma }}
T^{-1}_{\beta }T_{\beta }q_{a,\beta })\tilde{*} g_{\alpha } \}
$$
$$ = f_{\beta }\tilde{*} \{ \lim_a (T_{\beta } {\hat U}_{c_{\gamma }}
q_{a,\beta })\tilde{*}  g_{\alpha } \} \mbox{ and hence}$$
$$ \| ({\hat K}_{c_{\gamma }}f_{\beta
})\tilde{*} g_{\alpha } \|  \le $$
$$\overline{\lim}_a \| f_{\beta }\tilde{*} ([T_{\beta } {\hat U}_{c_{\gamma }}
q_{a,\beta }]\tilde{*} g_{\alpha }) \| \le \| f_{\beta } \| \|
T_{\beta } \| \| g_{\alpha } \| ~ \overline{\lim}_a \| [{\hat
U}_{c_{\gamma }} q_{a,\beta }]\tilde{*} \| \le \| f_{\beta } \| \|
g_{\alpha } \|
$$ for each $f, g \in {\cal E}$, since $ \| T \| \le 1$.
On the other hand, ${\hat K}_{c_{\gamma }^{-1}} =({\hat
K}_{c_{\gamma }})^{-1}$. Thus the inequalities $ \| {\hat
K}_{c_{\gamma }} \| \le 1$ and $ \| ({\hat K}_{c_{\gamma }})^{-1} \|
\le 1$ are satisfied for each $\gamma \in \Lambda $ and $c\in G$,
consequently, ${\hat K}_c$ is the isometry for each $c\in G$.
\par Applying Theorem 10 we get the statement of this lemma.
\par {\bf 12. Lemma.} {\it The mappings $(G,\tau _G^b)\ni c \to p_c\in (B^{\Lambda },\tau _B^b)$
for each $\beta $ and $(G,\tau _G^b)\ni c\mapsto \xi (c)\in (H,\tau
_H^b)$ of Lemma 11 are continuous homomorphisms, where $B= \{ x\in
{\bf F}: ~ |x|= 1 \} $ is the multiplicative group, the product
$B^{\Lambda }$ is in the box topology $\tau _B^b$, where $\tau ^b_G$
denotes the box topology on $G$ (see \S 9 \cite{bulopnlcga:12}).}
\par {\bf Proof.} These mappings are homomorphisms, since
$$p_{ch, \gamma } {\hat V}_{\xi _{\gamma }(c_{\gamma }h_{\gamma
})} = T_{\beta }{\hat U}_{c_{\gamma }h_{\gamma }}T^{-1}_{\beta } =
T_{\beta }{\hat U}_{c_{\gamma }}T^{-1}_{\beta }T_{\beta }{\hat
U}_{h_{\gamma }}T^{-1}_{\beta } = p_{c, \gamma }{\hat V}_{\xi
_{\gamma }(c_{\gamma })} p_{h, \gamma }{\hat V}_{\xi _{\gamma
}(h_{\gamma })}$$ for each $c, h \in G$, $~\beta \in \Lambda $ with
$\gamma = \phi (\beta )$, where $\xi (c) =\{ \xi_{\alpha }(c_{\alpha
}): ~ \alpha \in \Lambda \} $, $ ~ \xi _{\alpha } : G_{\alpha }\to
H_{\alpha }$ for each $\alpha \in \Lambda $. The mapping $\xi $ is
bijective, since for $\xi (c)=e_H\in H$, where $e_H$ is the neutral
element in $H$, one gets $p_{c,\gamma } I_{\cal F} = T_{\beta }
{\hat U}_{c_{\gamma }} T^{-1}_{\beta } $ and hence ${\hat
U}_{c_{\gamma }} =p_{c,\gamma }I_{\cal E}$, where $I_{\cal E}$
denotes the unit operator on $\cal E$. Therefore, $c=e_G$ and hence
$p_{c,\gamma }=1$ for each $\gamma $.
\par Then the mapping $G\ni c\mapsto {\hat U}_c$ is continuous
from $G$ in the box topology $\tau _G^b$ and relative to the strong
operator topology according to Proposition 10 \cite{bulopnlcga:12},
consequently, the mapping $H\ni t\mapsto {\hat V}_t$ is also
continuous, since $T$ and $T^{-1}$ are bounded linear operators.
\par Then for each $\epsilon = ( \epsilon_{\alpha } >0:
\alpha \in \Lambda )$ there exists a neighborhood $Y= \prod_{\alpha
\in \Lambda } Y_{\alpha }$ of $e_H$ in $(H, \tau _H^b)$ such that
each $Y_{\alpha }$ is an (open) neighborhood of the neutral element
$e_{\alpha }$ in $H_{\alpha }$ for which $\epsilon _{\alpha }/2 <
\lambda _{\alpha }( Y_{\alpha }) <\epsilon _{\alpha }$ for each
$\alpha \in \Lambda $, since $\lambda _{\alpha }$ is the
quasi-invariant borelian measure on $H_{\alpha }$ relative to the
dense subgroup $H_{\beta }$ and hence non-atomic. Moreover, if $Z$
is an arbitrary neighborhood of $e_H$ in $(H, \tau _H^b)$, then
there exists $Y$ such that $YY^{-1}\subseteq Z$. Then the function
$g = (g_{\alpha } = \chi _{Y_{\alpha }}: \alpha \in \Lambda )$
belongs to $\cal F$, where $\chi _{A_{\alpha }}$ denotes the
characteristic function of a subset $A_{\alpha }$ in $H_{\alpha }$.
Suppose that $p$ is a marked element in $B^{\Lambda }$. Let $t\in H$
be such that
$$(1)\quad \| p_{\beta } g_{\beta }\tilde{*}({\hat V}_{t_{\beta }}^*
g_{\alpha }) - g_{\beta }\tilde{*} g_{\alpha } \| < [\lambda _{\beta
}|_{Y_{\beta }}\tilde{*} \lambda _{\alpha }](Y_{\alpha })\mbox{,
 where}$$ $$ [\lambda _{\beta }|_{Y_{\beta }}\tilde{*} \lambda _{\alpha }](Y_{\alpha }):=
\int_{Y_{\beta }} \int_{Y_{\alpha }} \lambda_{\beta }(dx_{\beta
})\lambda_{\alpha }(\theta ^{\beta }_{\alpha }(x_{\beta })dx_{\alpha
}) ,$$ where $\theta ^{\beta }_{\alpha }: H_{\beta }\hookrightarrow
H_{\alpha }$ are embeddings (see \S 1). If $t_{\beta }\notin
Z_{\beta }$, then $s_{\beta }Y_{\beta }$ and $s_{\beta }t_{\beta
}Y_{\beta }$ are the disjoint subsets in the group $H_{\beta }$ for
each element $s_{\beta }$ in $H_{\beta }$, consequently,
$$\| p_{\beta } g_{\beta }\tilde{*}[{\hat V}_{t_{\beta }}^* g_{\alpha }] -
g_{\beta }\tilde{*} g_{\alpha } \| = \sup_{s_{\beta }\in H_{\beta }}
\int_{H_{\alpha }} | p_{\beta } [{\hat V}_{s_{\beta }t_{\beta
}}g_{\beta }] \tilde{*} g_{\alpha } (x_{\alpha }) - [{\hat
V}_{s_{\beta }}g_{\beta }] \tilde{*} g_{\alpha }(x_{\alpha }) |
\lambda _{\alpha }(dx_{\alpha })$$
$$=\sup_{s_{\beta }\in H_{\beta }} \int_{H_{\beta }} \int_{H_{\alpha }} | p_{\beta }
g_{\beta }(s_{\beta }t_{\beta }x_{\beta }) g_{\alpha } (\theta
^{\beta }_{\alpha }(x_{\beta })x_{\alpha }) | \lambda _{\beta
}(dx_{\beta }) \lambda _{\alpha }(dx_{\alpha })$$ $$ +
\sup_{s_{\beta }\in H_{\beta }} \int_{H_{\beta }} \int_{H_{\alpha
}}|g_{\beta }(s_{\beta }x_{\beta }) g_{\alpha }(\theta ^{\beta
}_{\alpha }(x_{\beta })x_{\alpha }) | \lambda _{\beta }(dx_{\beta })
\lambda _{\alpha }(dx_{\alpha })   \ge [\lambda _{\beta }|_{Y_{\beta
}}\tilde{*} \lambda _{\alpha }](Y_{\alpha })  .$$ Thus Inequality
$(1)$ implies that $t_{\beta }\in Z_{\beta }$. Hence the mapping $p
{\hat V}_{\xi _{\beta }(c_{\beta })}\mapsto \xi _{\beta }(c_{\beta
})=t_{\beta }\in H_{\beta }$, with $H_{\beta }$ in the topology
$\tau _{\beta }$, is continuous for each $\beta $, when linear
operators $p{\hat V}$ are considered relative to the strong operator
topology, since the set of all ($\mu _{\alpha }$-measurable) simple
functions is dense in $L^1_{G_{\beta }}(G_{\alpha })$. The mapping
$c_{\beta }\mapsto \xi _{\beta }(c_{\beta })$ is the composition of
three mappings $c_{\beta }\mapsto {\hat U}_{c_{\beta }}\mapsto
T_{\alpha }{\hat U}_{c_{\beta }}T^{-1}_{\alpha }=p_{c, \beta } {\hat
V}_{\xi _{\beta }(c_{\beta })} \mapsto \xi _{\beta }(c_{\beta
})=t_{\beta }$ which are continuous for each $\beta \in \Lambda $ as
it was proved above, consequently, the mapping $\xi : (G,\tau
_G^b)\to (H,\tau _H^b)$ is also continuous.
\par The mapping $c\mapsto p_c$ is continuous, since $c\mapsto p_cI$
is continuous as the composition of two uniformly bounded and
continuous mappings $T{\hat U}_cT^{-1}$ and ${\hat K}_{\xi (c)}$.
\par {\bf 13. Lemma.} {\it The mapping $\xi : G\to H$ is the
homeomorphism of $(G,\tau _G^b)$ onto $(H, \tau _H^b)$.}
\par {\bf Proof.} If $\{ \xi _{\beta } (x_{\beta ,b}) : b \} $ is a net
converging to $y_{\beta }\in H_{\beta }$, where $x_{\beta , b} \in
G_{\beta }$, then $\{ {\hat V}_{\xi _{\beta } (x_{\beta ,b})} : b
\}$ converges to ${\hat V}_{y_{\beta }}$ in the strong operator
topology. Therefore, $\{ T_{\alpha }^{-1}{\hat V}_{\xi _{\beta
}(x_{\beta ,b})}T_{\alpha } : b \}$ converges to $T_{\alpha
}^{-1}{\hat V}_{y_{\beta ,b}}T_{\alpha }$. From Lemma 11 we have the
equality
$$T_{\alpha }^{-1}{\hat V}_{\xi _{\beta }(x_{\beta ,b})}T_{\alpha } = p_{x_b
,\beta }^{-1} {\hat U}_{x_{\beta ,b}},$$ hence the net of operators
$\{ p_{x_b ,\beta }^{-1} {\hat U}_{x_{\beta ,b}}: b \} $ strongly
converges to $p_{\beta } {\hat U}_{x_{\beta }}$ for some $p_{\beta
}\in B$ and $x_{\beta }\in G_{\beta }$. Thus the equality $$p_{\beta
} T_{\alpha } {\hat U}_{x_{\beta }} T^{-1}_{\alpha } = {\hat
V}_{y_{\beta }}$$ is fulfilled with $y_{\beta }=\xi _{\beta
}(x_{\beta })$ and $p_{\beta }=p_{x,\beta }^{-1}$ for each $\beta
\in \Lambda $. This implies that $\xi _{\beta }(G_{\beta })$ is
closed in $H_{\beta }$ for each $\beta $ and hence $\xi (G)$ is
closed in $(H,\tau _H^b)$. \par The inverse operator $T^{-1}$ is
bounded (see \S 11). Then $T_{\alpha }^{-1}{\hat V}_{y_{\beta }}
T_{\alpha } = (sT_{\alpha })^{-1} {\hat V}_{y_{\beta }} (sT_{\alpha
})$ for each $s\in {\bf F}\setminus \{ 0 \} $. Hence without loss of
generality we can consider that $0 <\| T^{-1}_{\alpha } \| \le 1$
for each $\alpha \in \Lambda $. On the other hand, from the equality
$T^{-1}_{\alpha }{\hat V}_{y_{\beta }}T_{\alpha } = p_{x,\beta
}^{-1} {\hat U}_{x_{\beta }}$ with $x_{\beta }=\xi _{\beta
}^{-1}(y_{\beta })$ analogously to $\xi $ in \S 12 the continuity of
$\xi _{\beta }^{-1}: \xi _{\beta }(G_{\beta })\to G_{\beta }$
follows. \par Applying Lemmas 11, 12 and the proof in this section
above to $T^{-1}: {\cal F}\to {\cal E}$ we get that there exists a
continuous bijective homomorphism $\eta : (H,\tau _H^b)\to (G,\tau
_G ^b)$ such that $\eta (H)$ is closed in $(G,\tau _G ^b)$ and
\par $(1)$ ${\hat Q}_y=r_y{\hat U}_t$ for $t=\eta (y)$ and \par $(2)$ $r_y = \{
r_{y,\alpha }: |r_{y,\alpha }|=1 ~ \forall \alpha \in \Lambda \} \in
{\bf F}^{\Lambda }$, where ${\hat Q}_y = T^{-1} {\hat V}_yT$ for
each $y\in H$, $ ~ r: (G,\tau _G^b)\to B^{\Lambda }$ is a continuous
homomorphism. The operators ${\hat K}_c$ and ${\hat Q}_y$ are the
left meta-centralizers on ${\cal F}$ and ${\cal E}$ respectively for
each $c\in G$ and $y\in H$. But from 11$(1,2)$ it follows that $\eta
=\xi ^{-1}$ and $p_{\eta (y)}=r_y^{-1}$ for each $y\in H$, since
$\eta $ and $\xi $ are bijective homomorphisms. Therefore, Formulas
$(1,2)$ and 11$(1,2)$ imply that $\eta (\xi (G))=G$ and hence $\xi
(G)=H$.
\par {\bf 14. Theorem.} {\it Let $T: {\cal E}\to {\cal F}$ be an isomorphism of
normed algebras such that $Tf=(T_{\alpha }f_{\alpha }: \alpha \in
\Lambda )$, $T_{\alpha }: L^1_{G_{\beta }}(G_{\alpha }, \mu _{\alpha
},{\bf F})\to L^1_{H_{\beta }}(H_{\alpha }, \lambda _{\alpha },{\bf
F})$ and $\| T_{\alpha } \| \le 1$ for each $\alpha $, where ${\cal
F} = L^{\infty }(L^1_{H_{\beta }}(H_{\alpha }, \lambda _{\alpha
},{\bf F}): \alpha <\beta \in \Lambda )$ (see \S \S 11, 12). Then a
homeomorphism $\xi $ of topological groups exists from $(G,\tau
_G^b)$ onto $(H,\tau _H^b)$ and a continuous homomorphism $\psi :
G\to B^{\Lambda }$ such that
\par $(1)$ $T{\hat U}_xT^{-1}=\psi (x^{-1}) {\hat V}_{\xi (x)}$ and
\par $(2)$ $(Tf)_{\alpha }(\xi (x)) = \psi _{\beta }(x_{\beta }) f_{\alpha }(x_{\alpha })$
for each $x\in G$, $f\in {\cal E}$ and $\alpha \in \Lambda $ with
$\beta = \phi (\alpha )$, where $\psi (x)=(\psi _{\alpha }(x_{\alpha
}): \alpha \in \Lambda )$, $ ~ \psi _{\alpha }: G_{\alpha }\to B$,
$$~T_{\alpha }{\hat U}_{x_{\beta }}T^{-1}_{\alpha }=\psi _{\beta
}(x^{-1}_{\beta }) {\hat V}_{\xi _{\beta }(x_{\beta })}.$$ Moreover,
$T$ is an isometry.}
\par {\bf Proof.} We define a homomorphism $\psi (x)=p^{-1}_x$,
hence $\psi (x) = (\psi_{\alpha }(x_{\alpha }) = p^{-1}_{x,\alpha }
: \alpha \in \Lambda \} \in B^{\Lambda }$, hence $\psi _{\alpha } :
G_{\alpha }\to B$ is a character for each $\alpha \in \Lambda $.
From Lemmas 11-13 Statement $(1)$ of this theorem follows such that
$\xi : (G,\tau _G^b)\to (H,\tau _H^b)$ and $\xi ^{-1}: (H,\tau
_H^b)\to (G,\tau _G^b)$ and $\psi : G\to B^{\Lambda }$ are
continuous homomorphisms with $\xi (G)=H$.
\par If $S: {\cal E}\to {\cal F}$ is an isomorphism of normed
algebras such that $Sf=(S_{\alpha }f_{\alpha }: \alpha \in \Lambda
)$, $S_{\alpha }: L^1_{G_{\beta }}(G_{\alpha }, \mu _{\alpha },{\bf
F})\to L^1_{H_{\beta }}(H_{\alpha }, \lambda _{\alpha },{\bf F})$
and $\| S_{\alpha } \| \le 1$ for each $\alpha $ such that $S$
satisfies Equality $(2)$: \par $(Sf)_{\alpha }(\xi (x)) = \psi
_{\beta }(x_{\beta }) f_{\alpha }(x_{\alpha })$ for each $x\in G$
and $f\in {\cal E}$, then $(S^{-1} g)_{\alpha } (x) = \psi _{\beta
}(x^{-1}_{\beta }) g_{\alpha }(\xi _{\alpha }(x_{\alpha }))$ for
each $g\in {\cal F}$ and $x\in G$. Therefore one infers that
$$(S_{\alpha }{\hat U}_{c_{\beta }}S^{-1}_{\alpha }g_{\alpha })(\xi
_{\alpha }(x_{\alpha })) = \psi _{\beta }(x_{\beta }) ({\hat
U}_{c_{\beta }}S^{-1}_{\alpha }g_{\alpha })(x_{\alpha }) =$$ $$ \psi
_{\beta }(x_{\beta }) (S^{-1}_{\alpha }g_{\alpha })(\theta ^{\beta
}_{\alpha }(c_{\beta }) x_{\alpha }) = \psi _{\beta }(x_{\beta })
\psi _{\beta }(x_{\beta }^{-1} c_{\beta }^{-1}) g_{\alpha }(\theta
^{\beta }_{\alpha }(\xi _{\beta }(c_{\beta }))\xi _{\alpha
}(x_{\alpha }) )
$$
$$ =\psi _{\beta }(c_{\beta }^{-1}) g_{\alpha }(\theta ^{\beta }_{\alpha
}(\xi _{\beta }(c_{\beta }))\xi _{\alpha }(x_{\alpha }) ) = \psi
_{\beta }(c_{\beta }^{-1}) ({\hat U}_{\xi _{\beta }(c_{\beta
})}g_{\alpha })(\xi _{\alpha }(x_{\alpha }) ),$$ consequently,
$S_{\alpha } {\hat U}_{c_{\beta }} S^{-1}_{\alpha } = \psi _{\beta
}(c_{\beta }^{-1}) {\hat U}_{\xi _{\beta }(c_{\beta })} $ for each
$c\in G$, $\alpha \in \Lambda $ with $\beta = \phi (\alpha )$, where
embeddings $H_{\beta }\hookrightarrow H_{\alpha }$ also are denoted
by $\theta ^{\beta }_{\alpha }$ for the notation simplicity (see \S
1). This means that $S{\hat U}_cS^{-1} =T{\hat U}_cT^{-1}$ and hence
$$(3)\quad T^{-1}_{\alpha }S_{\alpha }{\hat U}_{c_{\beta }}={\hat
U}_{c_{\beta }}T^{-1}_{\alpha }S_{\alpha }$$ for each $\alpha \in
\Lambda $ with $\beta = \phi (\alpha )$. In view of Lemmas 11-13 and
the conditions of this theorem the linear operators $T$, $T^{-1}$,
$S$ and $S^{-1}$ are continuous. Thus the operator
\par $(4)$ $T^{-1}S=:Y$ \\ is the isomorphism of the algebra
$\cal E$ onto itself commuting with all operators ${\hat U}_c$ such
that $Y$ and $Y^{-1}$ are continuous. As in \S 13 it is sufficient
to consider the case $ 0 < \| Y_{\alpha } \| \le 1$ for each $\alpha
\in \Lambda $, since ${\hat U}_{c_{\beta }} = Y_{\alpha }^{-1}{\hat
U}_{c_{\beta }} Y_{\alpha } = (kY_{\alpha })^{-1} {\hat U}_{c_{\beta
}} (k Y_{\alpha } ) $ for every $k\in {\bf F}\setminus \{ 0 \} $,
$\alpha \in \Lambda $ with $\beta = \phi (\alpha )$ and $c\in G$.
Take $f, q\in {\cal E}$ and consider the left meta-centralizer $A$
defined by a radonian measure $\nu _{\alpha }\in M_t(G_{\alpha
},{\bf F})$ such that
\par $(5)$ $\nu _{\alpha }(dx_{\alpha })=q_{\alpha }(x_{\alpha })\mu
_{\alpha }(dx_{\alpha })$ \\ for each $\alpha \in \Lambda $, that is
$Af=\nu \tilde{\star }f$. On the other hand, $$(6)\quad (Af)_{\alpha
}(x_{\alpha }) = \int_{G_{\beta }} q_{\beta }(y_{\beta })[{\hat
U}_{y_{\beta }}f_{\alpha }(x_{\alpha })]\mu _{\beta }(dy_{\beta
}),$$ that is relative to the strong operator topology
$$(7)\quad A_{\alpha } =\int_{G_{\beta }} q_{\beta }(y_{\beta })
{\hat U}_{y_{\beta }} \mu _{\beta }(dy_{\beta })$$ for each $\alpha
\in \Lambda $ with $\beta =\phi (\alpha )$, where $Af =(A_{\alpha
}f_{\alpha }: \alpha \in \Lambda )$. In each Banach space
$L^1_{G_{\gamma }}(G_{\beta },\mu _{\beta },{\bf F})$ the space of
($\mu _{\beta }$-measurable) simple functions $\sum_{j=1}^n v_j \chi
_{Z_j}$ is dense, where $v_j\in \bf F$ is a constant and $Z_j$ is a
$\mu _{\beta }$-measurable subset in $G_{\beta }$ for each
$j=1,..,n$, $ ~ n\in \bf N$. Therefore, from Formulas $(3-7)$ it
follows that $$YAf= Y(q\tilde{\star }f) = (Yq) \tilde{\star }(Yf) =
AYf =q \tilde{\star }(Yf),$$ consequently, $Yq=q$ for each $q\in
\cal E$, since $f\in {\cal E}$ is arbitrary. Thus $Y=I_{\cal E}$ and
hence $T=S$, where $I_{\cal E}$ denotes the unit operator on ${\cal
E}$. From this Formula $(2)$ follows. The last statement follows
from Formulas $(2)$ and 3$(1)$.

\par {\bf 15. Remark.} The results of this paper can be used for
further studies of non locally compact group algebras,
representations of groups, completions and extensions of groups,
etc.

\end{document}